\documentclass[reqno,11pt]{amsart}

\usepackage{color}
 \usepackage{mathtools}

\usepackage{graphicx}

 \usepackage{amsmath}  
 \usepackage{amssymb}
\usepackage{bm}

\newtheorem{theorem}{Theorem} 
\newtheorem{lemma}[theorem]{ Lemma}
\newtheorem{corollary}[theorem]{ Corollary}

\theoremstyle{definition}

\theoremstyle{remark}
\newtheorem{remark}[theorem]{Remark}

\makeatletter 
 \def\dashint{\operatorname{\,\,\,\mathclap{\int} \kern-.23em\text{\bf--}\!\!}}
 \makeatother

\makeatletter 

\def\dashnorm{\,\,\text{\bf--}\kern-.5em\|}
\def\ninf{\qopname\relax\@empty{inf\phantom{p}\!\!\!}}

 \makeatother

\newcommand\bB{\mathbb{B}}

\newcommand\bM{\mathbb{M}}
\newcommand\bR{\mathbb{R}}

\begin{document}

\title[A remark on a paper of
F. Chiarenza and M. Frasca]{A remark on a paper of
F. Chiarenza and M. Frasca}

\author{N.V. Krylov}
 
\email{nkrylov@umn.edu}
\address{127 Vincent Hall, University of Minnesota,
 Minneapolis, MN, 55455}



\begin{abstract}
In 1990 F. Chiarenza and M. Frasca published
a paper in which they generalized a result
of C. Fefferman on estimates of
the integral of $|bu|^{p}$ through the integral
of $|Du|^{p}$ for   $p>1$.
Formally their proof is valid only for $d\geq 3$.
We present here further generalization
with  a different proof in which
  $D $ is replaced with the fractional
power   of the Laplacian for any dimension
$d\geq 1$.
\end{abstract}

\maketitle 

Let an integer $d\geq 1$ and let $\bR^{d}$
be a Euclidean space of points $x=(x^{1},...,x^{d})$.
Fix
$\alpha\in(0,d)$ and consider the Riesz
potential
$$
R_{\alpha}f(x)=\int_{\bR^{d}}\frac{f(x+y)}
{|y|^{d-\alpha}}\,dy.
$$

We denote by $B_{r}(x)$ the open ball of radius
$r$ centered at $x$, $B_{r}=B_{r}(0)$, $\bB_{r}$
the collection of $B_{r}(x)$, $S_{1}=\{|x|=1\}$.
Our main result is the following,
in which $r,p,A$ are some numbers and $b=b(x)$
is a measurable function.

\begin{theorem} 
                            \label{theorem 10.7.1}

Assume $\alpha\leq r$, $1<r<p\leq d$, $b\geq0$, $f\in L_{r}$, and
for {any\/} $\rho>0$ and $B\in\bB_{\rho}$.
$$
\Big(\dashint_{B}b^{p}\,dx\Big)^{1/p}\leq A\rho^{-\alpha}. 
$$
Then
\begin{equation}
                                \label{9.25.1}
I:=\int_{\bR^{d}}b^{r}|R_{\alpha}f|^{r}\,dx\le N(\alpha,d,r,p)A^{r}
 \int_{\bR^{d}} |f|^{r}\,dx.
\end{equation}
\end{theorem}

Below by $N$ we denote generic constants
depending only on $\alpha,d,r,p,q$.

\begin{corollary}
                       \label{corollary 10.7.1}
If $u\in C_{0}^{\infty}(\bR^{d})$, then
$$
\int_{\bR^{d}}b^{r}|u|^{r}\,dx\le  NA^{r}
 \int_{\bR^{d}} \big|(-\Delta)^{\alpha/2}u|^{r}\,dx
$$
and, if $\alpha=1$ and hence $d\geq2$, (the Chiarenza-Frasca result)
$$
\int_{\bR^{d}}b^{r}|u|^{r}\,dx\le  NA^{r}
 \int_{\bR^{d}}  |Du|^{r}\,dx.
$$
\end{corollary}
Indeed, $f:=(-\Delta)^{\alpha/2}u$ satisfies
$f\in L_{r}$ and $R_{\alpha}f=u$ and the $L_{r}$-norms of $Du$ and $(-\Delta)^{1/2}u$ are equivalent.

\begin{remark}
The author have used the Chiarenza-Frasca theorem
in a few papers leading to \cite{11}
about strong solutions of It\^o's equations.
Theorem \ref{theorem 10.7.1} paves the way to
treat equations driven by L\'evy rather than Wiener processes.
\end{remark}
 
We prove Theorem \ref{theorem 10.7.1}
adapting to the ``elliptic'' setting
the proof of Theorem 4.1 of \cite{Kr_23}.
We need two auxiliary and certainly well-known
results in which $\bM$ is the Hardy-Littlewood
maximal operator.

\begin{lemma}
                      \label{lemma 9.23.1}
If   $1< q\leq p$, it holds that
$$
R_{\alpha}(b^{q})\leq NA (\bM (b^{q}))^{1-1/ q }.
$$

\end{lemma}

Proof. We have
$$
R_{\alpha}b^{q}(0)=N\int_{0}^{\infty}
r^{\alpha-1}\int_{S_{1}}b^{q}(r\theta)\,\sigma(d\theta)\,dr
$$
$$
=N\int_{0}^{\infty}r^{\alpha-d}
\frac{d}{dr}\int_{B_{r}}b^{q}\,dx\,dr 
$$
$$
\leq N\int_{0}^{\infty}r^{\alpha-d-1}
\int_{B_{r}}b^{q}\,dx\,dr
=N\int_{0}^{\rho}+N\int_{\rho}^{\infty}
$$
$$
\leq N\rho^{\alpha}\bM (b^{q})+N\rho^{\alpha-q\alpha}A^{q},
$$
where we used that
$$
\dashint_{B_{r}}b^{q}\,dx\leq
\Big(\dashint_{B_{r}}b^{p}\,dx\Big)^{q/p}\leq A^{q}r^{-q\alpha}.
$$
 For
$$
\rho^{-q\alpha}=\bM (b^{q})/A^{-q},\quad \rho^{\alpha}
=A[\bM (b^{q})]^{-1/q }
$$
we get the result. \qed

\begin{lemma}
                        \label{lemma 9.25.1}
For any $\rho>0$
\begin{equation}
                         \label{9.25.2}
I:=\int_{\bR^{d }} b ^{p}\bM I_{B_{\rho}}\,dx\leq NA^{p}\rho^{d-p\alpha}.
\end{equation}

\end{lemma}

Proof. We have
$$
\bM I_{B_{\rho}}\leq N(I_{B_{ \rho}}+
I_{|x|> \rho}\frac{\rho^{d}}{|x|^{d}})
$$
and
$$
\int_{\bR^{d }} b ^{p}I_{B_{ \rho}}\leq NA^{p}\rho^{d-p\alpha},
$$
$$
\rho^{d}\int_{\rho}^{\infty}r^{-d}
\frac{d}{dr}\int_{B_{r}}b^{p}\,dx\,dr
\leq d\rho^{d}\int_{\rho}^{\infty}
r^{-d-1}\int_{B_{r}}b^{p}\,dx\,dr
$$
$$
\leq N\rho^{d}A^{p}\int_{\rho}^{\infty}r^{-1
-p\alpha}\,dr=N\rho^{d-p\alpha}A^{p}.
$$
This yields the result. \qed

{\bf Proof of Theorem \ref{theorem 10.7.1}}. It suffices to concentrate on $f\geq0$. Then first assume that $b^{p}\in A_{1}$,
that is $\bM (b^{p})\leq Nb^{p}$.
Observe that for $v=R_{\alpha}f$ we have 
\begin{equation}
                           \label{5.30.2}
I= \int_{\bR^{d }}\big(b^{r}v^{r-1}\big) R_{\alpha}f
\,dx=\int_{\bR^{d}}R_{\alpha}\big(b^{r}v^{r-1}\big)  f
\,dx\leq \|f\|_{L_{r}}\big\|R_{\alpha}\big(b^{r}v^{r-1}\big) \big\|_{L_{r'}},
\end{equation}
where $r'=r/(r-1)$.

Next, take $\gamma>0$, such that $(1+\gamma)r
\leq p$, $1+\gamma r'\leq p$,
 and $r\geq 1+\gamma$. Note that
$$
R_{\alpha}\big(b^{r}v^{r-1}\big)=
R_{\alpha}\big(b^{1+\gamma}\big(b^{r-1-\gamma}
v^{r-1}\big)\big)
$$
$$
\leq\Big(R_{\alpha}\big(b^{(1+\gamma)r})\Big)^{1/r}
\Big( R_{\alpha}\big(b^{r-\gamma r'}v^{r})\Big)^{(r-1)/r}.
$$
It follows that  
$$
\big\|
R_{\alpha}\big(b^{r}v^{r-1}\big)\big\|_{L_{r'}}
\leq\Big(\int_{\bR^{d}}b^{r-\gamma r'}v^{r}
R_{\alpha}\Big[\Big(R_{\alpha}\big(b^{(1+\gamma)r})\Big)^{1/(r-1)}\Big]\,dx\Big)^{(r-1)/r}.
$$
Now in light of \eqref{5.30.2} we see that, to prove the theorem in our particular case, it only remains to show that
\begin{equation}
                             \label{5.30.4}
R_{\alpha}\Big[\Big(R_{\alpha}\big(b^{(1+\gamma)r})\Big)^{1/(r-1)}\Big]\leq Nb^{\gamma r'}A^{r'} .
\end{equation}

By observing that
$1<(1+\gamma)r\leq p$ 
and using Lemma \ref{lemma 9.23.1} we get  that 
$$
R_{\alpha}\big(b^{(1+\gamma)r}\big)\leq NA \big(  \bM\big(b^{(1+\gamma)r}\big)\big)^{1-1/(r+\gamma r )},
$$
where by assumption 
and H\"older's inequality
$$
\big(  \bM\big(b^{(1+\gamma)r}\big)\big)^{1-1/(r +\gamma r )}=\big[\big(  \bM\big(b^{(1+\gamma)r}\big)\big)^{1/(r+\gamma r)}\big]^{
(1+\gamma)r-1}
$$
$$
\leq Nb^{(1+\gamma)r-1}=Nb^{r-1+\gamma r}.
$$
  Hence,
$$
R_{\alpha}\Big[\Big(R_{\alpha}\big(b^{(1+\gamma)r})\Big)^{1/(r-1)}\Big]\leq
NA^{1/(r - 1)}
R_{\alpha}b^{1+\gamma r'}.
$$
By Lemma \ref{lemma 9.23.1}
$$
R_{\alpha}b^{1+\gamma r'}\leq NA (\bM(b^{1+\gamma r'}))
^{1-1/(1+\gamma r')}\leq NAb^{\gamma r'}.
$$
This yields \eqref{5.30.4} and proves the 
lemma in our particular case.

We now get rid of the assumption that 
$\bM( b ^{p})\leq N  b ^{p}$ as in \cite{CF_90}. 
For $p_{0}=(r+p)/2$, $p_{1}=(r+p_{0})/2$ we have $b^{p_{1}}\leq 
(\bM(b^{p_{0}}))^{p_{1}/p_{0}}:=\tilde b^{p_{1}}$ and since $p_{1}/p_{0}<1$, $\tilde b^{ p_{1}}$ is an $A_{1}$-weight
with  the $A_{1}$-constant depending only on
$p_{1}/p_{0}$  (see, for instance, \cite{GR_85}, p.~158).
Therefore, \eqref{9.25.1} holds with $\tilde b$
in place of $b$ and it only remains to show that    for any $x,\rho$,
\begin{equation}
                             \label{1.19.1}
\int_{B_{\rho}( x)}\tilde b^{p_{1}}\,dx
\leq  N\rho^{d-p_{1}\alpha}A^{p_{1}}.
\end{equation}
Of course, we may assume that $x=0$.

Then by H\"older's inequality we see that
the left-hand side of \eqref{1.19.1} is less than 
$$
N\rho^{d(p-p_{1})/p}\Big(\int_{\bR^{d }}
(\bM(b^{p_{0}}))^{p/p_{0}}I_{B_{\rho}}\,dx\Big)^{p_{1}/p},
$$
where the integral 
by a Fefferman-Stein Lemma 1, p.~111 of
\cite{FS_71} and the fact that
$p/p_{0}>1$ is dominated by
$$
N\int_{\bR^{d }} b ^{p}\bM I_{B_{\rho}}\,dx
\leq NA^{p}\rho^{d-p\alpha},
$$
where we used Lemma \ref{lemma 9.25.1}.
Hence, 
$$
\int_{B_{\rho}}\tilde b^{p_{1}}\,dx
\leq N\rho^{d(p-p_{1})/p}A^{p_{1}}\rho^{p_{1}d/p-p_{1}\alpha}
$$
which is \eqref{1.19.1}. 
The theorem is proved.  \qed

\end{document}